\documentstyle[12pt]{article}
\newtheorem{theorem}{THEOREM}[section]
\newtheorem{definition}[theorem]{DEFINITION}

\newtheorem{remarks}[theorem]{REMARKS}

\newtheorem{cor}[theorem]{COROLLARY}

\title
{\normalsize\bf
\vskip 2truecm
GENERALIZED INVERSES AND THE MAXIMAL RADIUS OF REGULARITY 
OF A FREDHOLM OPERATOR
}
\author
{\normalsize
Catalin BADEA
\and
		{\normalsize Mostafa MBEKHTA}
}

\date{}
\begin{document}
\maketitle
\thispagestyle{empty}

\baselineskip=12pt
\begin{quote}
\noindent Operators possessing analytic generalized inverses 
satisfying the resolvent
identity are studied. Several characterizations and necessary conditions
are obtained. The maximal radius of regularity
for a Fredholm operator $T$ is computed in terms of the spectral radius of
a generalized inverse of $T$. This provides a partial answer to a conjecture 
of J. Zem\'anek.
\end{quote}

\vskip 1truecm
\baselineskip=15pt
	
\thispagestyle{empty}
\section{Introduction}
			Let $X$ be a Banach space. We will denote by $C(X)$ the set
of all closed operators with a dense domain and by $B(X)$ the algebra of bounded
operators from $X$ into itself. For an operator $T \in C(X)$,
we denote by $D(T), N(T)$ and $R(T)$ the domain, the kernel and, respectively,
the range of $T$. The identity operator will be denoted by $I$. 
We introduce the following definition for the generalized
inverse of $T$ (relative inverse or pseudo-inverse are names 
also used in the
literature).

\begin{definition}
		The operator $T \in C(X)$ possess a \it{generalized inverse\/} if there exists an
		operator $S \in B(X)$ such that $R(S) \subseteq D(T)$ and
				\begin{enumerate}
				\item $TST = T$ on $D(T)$.
				\item $STS = S$ on $X$.
				\item $ST$ is continuous \label{def11-3}.
				\end{enumerate}
\end{definition}

\noindent In this case we will say that $S$ is the generalized inverse of $T$.

The following are some simple remarks about generalized inverses.

\begin{remarks}  (cf. \cite{C})  \end{remarks}
\begin{enumerate}
				\item Using the closed-graph theorem, we get that the operator $TS$ is continuous.
				\item The operator $TS$ is a projection onto $R(T)$ such that
										$N(TS) = N(S)$ and $R(TS) = R(T)$.
				\item The condition~\ref{def11-3} in the above definition
										shows that the operator $ST$ can be extended to a bounded projection
										onto $\overline{R(S)}$ such that $N(ST) = N(T)$ and $R(ST) = \overline{R(S)}$.
				\item The operator $T$ possess a generalized inverse if and only if $N(T)$
										and $R(T)$ have topological complements in $X$.
				\end{enumerate}
In what follows we will denote $P = TS$ and $Q = ST$.

\smallskip

The following is an open problem : Suppose that $U \subset reg(T)$, that is, for any $\lambda \in U$, 
$T - \lambda I$ possesses an analytic 
generalized inverse in a suitable neighborhood $V_{\lambda}$ of $\lambda$, 
where $U$ is an open, connected subset of {\bf C}. Does 
a generalized resolvent of $T$ on $U$ always exist ? A generalized resolvent of $T$ on 
$U$ is an operator-valued
function $Rg(T,\lambda)$ on $U$ such that $Rg(T,\lambda)$ is a generalized 
inverse
of $T - \lambda I$ for all $\lambda \in U$ and $Rg(T,\lambda)$ satisfies 
the resolvent
identity. This open problem is mentioned in \cite{Saph}, \cite{Mbe2},
\cite{Mbe3}, \cite{Sch}, \cite{Mul}. Note also that the corresponding problem for 
left or right invertible operators is also 
open (cf. \cite{AC2}, \cite{Her}, \cite{Zem}).

The paper is organized as follows. In the next section we recall some known 
results and we prove 
several characterizations for generalized 
resolvents. These results will be used in the proof of the main results 
of this paper (section 3). A consequence of one of it gives a formula 
for the maximal radius of 
regularity of a 
Fredholm operator which answers, in a particular case, a question of 
Zem\'anek  \cite{Zem}.  


\section{Generalized inverses and generalized resolvents}
We start by recalling some notation and definitions.

\begin{definition}
		Let $U$ be a subset of {\bf C}. We will say that $T\in C(X)$ possess a
	 generalized inverse on $U$ if $T - \lambda I$ possess a generalized inverse
		for every $\lambda \in U$.
\end{definition}

We will denote by $reg(T)$ the set of all complex numbers $\lambda$ for which
$T$ possess an {\it analytic\/} generalized inverse in a neighborhood of 
$\lambda$. Then
$$\sigma_{g}(T) = {\bf C}\setminus reg(T)$$ will 
denote the {\it generalized spectrum\/} of $T$. Several properties of the 
classical
spectrum $\sigma (T)$ remain true in the case of 
the generalized one (cf. \cite{Mbe1, Mbe2}). 

\begin{definition}
		An operator $T \in C(X)$
is called regular if $T$ possess a generalized inverse 
and $N(T^{n}) \subset R(T)$,
for all $n \geq 0$.
\end{definition}

It is easy to see that the condition ``$N(T^{n}) \subset R(T)$,
for all $n \geq 0$'' is equivalent to the following one : ``$N(T) \subset R(T^{m})$,
for all $m \geq 0$''.

\smallskip
The following result gives a characterization of $reg(T)$ in terms of regular
operators.

\begin{theorem}[(cf. \cite{Mbe1})]
For an operator $T \in C(X)$, we have $\lambda _{0} \in reg(T)$
if and only if $T - \lambda_{0}I$ is regular.
\end{theorem}

In Section 3, we will need the following result for regular operators. 
It was proved in \cite{Mbe0} for bounded operators on a
Hilbert space but the proof there remains valid for this more general case.

\begin{theorem}[(cf. \cite{Mbe0})]
Let $T \in C(X)$ be a regular operator and $S \in B(X)$ 
a generalized inverse
of $T$. Then, for every $i \geq 1$, we have
$$T^{i}S^{i}T^{j} = 
\left\{ 
\begin{array} {r@{\quad:\quad}l} T^{i}S^{i-j} & 0 \leq j \leq i \\ 
T^{j} & i \leq j \end{array} \right. $$
\noindent and
$$T^{j}S^{i}T^{i} = 
\left\{ 
\begin{array} {r@{\quad:\quad}l} S^{i-j}T^{i} & 0 \leq j \leq i \\ 
T^{j} & i \leq j \end{array} \right. $$
In particular we have $T^{n}S^{n}T^{n} = T^{n}$ for all $n \geq 1$.
\end{theorem}

Let $U \subset {\bf C}$ be an open set.
We will say that the operator-valued function $f$ defined on 
$U$ {\it satisfies the resolvent identity\/} on $U$ if
$$f(\lambda) - f(\mu) = (\lambda - \mu)f(\lambda)f(\mu)$$
for all $\lambda$ and $\mu$ belonging to the {\it same \/} connected component of $U$.

\begin{definition}
Let $U \subset {\bf C}$ be an open set. The operator $T \in C(X)$ is said to 
possess a generalized resolvent on $U$ if $T$ possess a generalized
inverse on $U$ satisfying the resolvent identity in $U$.
\end{definition} 

\noindent Note that a generalized resolvent $Rg(T,\lambda)$
of $T$ in a connected $U$ is analytic on $U$. According to \cite[page 184]{HiPh}, a function satisfying
the resolvent identity on an open set is locally analytic.

\smallskip

The following result characterizes generalized resolvents (on connected,
open sets) among the analytic generalized inverses. It generalizes 
a result from \cite{Zem} and will be used in the 
next section in the proof of the main result.

\begin{theorem} \label{cha1}
Let $U$ be an open, connected subset of {\bf C}, $0 \in U$. Let $T \in C(X)$ 
possessing an analytic generalized inverse 
$R(\lambda) = \sum_{n=0}^{\infty}\lambda^{n}T_{n}$ on $U$. Denote
$P(\lambda) = (T - \lambda I)R(\lambda)$ and 
$Q(\lambda) = R(\lambda)(T - \lambda I)$ for $\lambda \in U$. 
The following conditions are equivalent :
	\begin{enumerate}
				\item[(i)] $T_{n} = T_{0}^{n+1}$ for all $n \geq 1$.
				\item[(ii)]$N(P(\lambda)) = N(TT_{0})$ and $R(Q(\lambda)) = R(T_{0}T)$
								for all $\lambda \in U$.
				\item[(iii)] There exist two closed subspaces $Z$ and $W$ of $X$ such that
						$N(P(\lambda)) = Z$ and $R(Q(\lambda)) = W$ for all $\lambda \in U$.
				\item[(iv)] $R(\lambda)$ is a generalized resolvent of $T$ on $U$.
				\item[(v)] $R(\lambda) - R(0) = \lambda R(\lambda)R(0)$, for all 
															$\lambda \in U$.
				\end{enumerate}
\end{theorem} 

{\it PROOF.\/} (i) $\Rightarrow$ (ii) : Since $R(\lambda)$ is a 
generalized inverse of $T$ on $U$ it follows, in particular, that
$T_{0}$ is a generalized inverse of $T$. Since $T_{n} = T_{0}^{n+1}$ for 
all $n \geq 1$, we have
$$Q(\lambda) = R(\lambda)(T - \lambda I) 
= T_{0}T - \sum_{n=1}^{\infty}\lambda^{n}T_{0}^{n}(I - T_{0}T)$$
and
$$P(\lambda) = (T - \lambda I)R(\lambda) 
= TT_{0} - (I - TT_{0})\sum_{n=1}^{\infty}\lambda^{n}T_{0}^{n}$$
for every $\lambda \in U$. We prove first the equality 
$N(P(\lambda)) = N(T_{0})$. By the Remark 1.2, (2), we will
have $N(P(\lambda)) = N(TT_{0})$. 
Let $u \in N(P(\lambda))$. Then
$$0 = P(\lambda)u = TT_{0}u - (I - TT_{0})\sum_{n=1}^{\infty}\lambda^{n}T_{0}^{n}u.$$
Therefore 
$$TT_{0}u = (I - TT_{0})\sum_{n=1}^{\infty}\lambda^{n}T_{0}^{n}u.$$
Applying $T_{0}$ to both sides we get $T_{0}u = 0$ (we have used the equality 
$T_{0}TT_{0} = T_{0}$). Therefore $N(P(\lambda)) \subseteq N(T_{0})$. 
Now let $u \in N(T_{0})$. For every $\lambda 
\in U$, we have 
$$P(\lambda)u = TT_{0}u - (I - TT_{0})\sum_{n=1}^{\infty}\lambda^{n}T_{0}^{n}u = 0.$$
Hence $N(P(\lambda) = N(TT_{0})$.

We prove now that $R(Q(\lambda)) = R(T_{0}T)$. Let $u = Q(\lambda)u \in 
R(Q(\lambda))$ for a fixed $\lambda \in U$. Then
$$u = T_{0}Tu - \sum_{n=1}^{\infty}\lambda^{n}T_{0}^{n}(I - T_{0}T)u$$
$$= T_{0}Tu - \sum_{n=1}^{\infty}\lambda^{n}
T_{0}^{n}(I - T_{0}T)\left[ T_{0}\left( Tu - \sum_{n=1}^{\infty}
\lambda^{n}T_{0}^{n-1}(I - T_{0}T) \right) u\right]$$
$$= T_{0}Tu.$$
\noindent Hence $R(Q(\lambda)) \subseteq R(T_{0}T)$. 

Conversely, if $u = T_{0}Tu \in R(T_{0}T)$, then 
$$Q(\lambda)u = T_{0}Tu - \sum_{n=1}^{\infty}\lambda^{n}T_{0}^{n}(I - T_{0}T)u$$
$$= u - \sum_{n=1}^{\infty}\lambda^{n}T_{0}^{n}(I - T_{0}T)T_{0}Tu 
= u \in R(Q(\lambda)).$$
Therefore $R(Q(\lambda)) = R(T_{0}T)$ for every $\lambda \in U$.

(ii) $\Rightarrow$ (iii) is clear.

(iii) $\Rightarrow$ (iv) : Let $\lambda , \mu \in U$. Then
$$(\lambda - \mu)R(\lambda)R(\mu) = R(\lambda)(\lambda - \mu)R(\mu)$$
$$= R(\lambda)\left[ (T - \mu I) - (T - \lambda I)\right] R(\mu)
 = R(\lambda)P(\mu) - Q(\lambda)R(\mu).$$
Suppose now that $N(P(\lambda)) = Z$ for all $\lambda \in U$. Then 
we have $R(\lambda)(I - P(\mu)) = 0$
and so $R(\lambda) = R(\lambda)P(\mu)$ for all $\lambda , \mu \in U$.

Supposing also that $R(Q(\lambda)) = W , \lambda \in U$, we obtain 
$(I - Q(\lambda))R(\mu) = 0$ for all $\lambda , \mu \in U$.
Therefore
$$(\lambda - \mu)R(\lambda)R(\mu) = R(\lambda)P(\mu) - Q(\lambda)R(\mu)
 = R(\lambda) - R(\mu).$$

(iv) $\Rightarrow$ (v) is clear.

(v) $\Rightarrow$ (i) : Suppose that 
$R(\lambda) - R(0) = \lambda R(\lambda)R(0)$. Then, for every $\lambda \in U$,
we have
$$\sum_{n=0}^{\infty}\lambda^{n}T_{n} - T_{0} 
= \lambda \sum_{n=0}^{\infty}\lambda^{n}T_{n} T_{0}.$$
Therefore
$$\sum_{n=1}^{\infty}\lambda^{n}T_{n} - \sum_{n=0}^{\infty}\lambda^{n+1}T_{n}T_{0} = 0$$
and thus
$$\sum_{n=0}^{\infty}\left( T_{n+1} - T_{n}T_{0}\right) \lambda^{n+1} = 0.$$
Hence $T_{n+1} - T_{n}T_{0} = 0$ for all $n \geq 0$, yielding 
$T_{n} = T_{0}^{n+1}$
for all $n \geq 0$. \hfill $\diamondsuit$

\smallskip
The next result is a characterization of generalized resolvents, without
assuming the power-series development around $0$. We will apply this 
criterion in Example 3.6.

\begin{theorem}   \label{cha2}
Let $T \in C(X)$ and let $U$ be an open, connected subset of {\bf C}. The following conditions
are equivalent :
	\begin{enumerate}
				\item[(i)] There exist two families of projections
															$P(\lambda)$ and $Q(\lambda)$, $\lambda \in U$, continuous in 
										$\lambda$, such that
	$$R(P(\lambda)) = R(T - \lambda I) , N(Q(\lambda)) = N(T - \lambda I) ; \lambda \in U$$
	and
	$$P(\lambda)P(\mu) = P(\lambda) , Q(\lambda)Q(\mu) = Q(\mu) ; \lambda , \mu \in U.$$
				\item[(ii)] There exists a generalized resolvent of $T$ on $U$.
				\item[(iii)] There exists an analytic generalized inverse $R(\lambda)$ of 
		$T$ in $U$, satisfying $R'(\lambda) = R(\lambda)^{2}$, for all 
															$\lambda \in U$.
				\end{enumerate}

\end{theorem}

{\it PROOF.\/} (i) $\Rightarrow$ (ii) : Let $u \in X$ and $\lambda \in U$.
Then $P(\lambda)u \in R(T - \lambda I)$. Therefore, there exists $v \in D(T)$
such that $P(\lambda)u = (T - \lambda I)v$. Set $Rg(T, \lambda)u = Q(\lambda)v$.
Firstly, we show that $Rg(T, \lambda)$ is well-defined. Indeed, if $w \in D(T)$ is such
that $(T - \lambda I)w = P(\lambda)u = (T - \lambda I)v$, then $$v - w \in N(T - \lambda I) 
= N(Q(\lambda)).$$ 
Therefore $Q(\lambda)(v - w) = 0$ and thus $Q(\lambda)v = Q(\lambda)w$.
Hence $Rg(T,\lambda)$ does not depend on the choice of $v$.

We show that $Rg(T,\lambda)$ is a generalized inverse of $T - \lambda I$. 
For all $u \in X$, we have $(I - Q(\lambda))u \in N(T - \lambda I) \subset D(T)$. Also 
$(T - \lambda I)(I - Q(\lambda)) = 0$ and $Q(\lambda)(D(T)) \subset D(T)$. 
We obtain $T - \lambda I = (T - \lambda I)Q(\lambda)$ and $R(Rg(T,\lambda)) \subset D(T)$.
Using the definition of $Rg(T,\lambda)$, we have
$$(T - \lambda I)Rg(T,\lambda)(T - \lambda I)u = (T - \lambda I)Q(\lambda)u \\
 = (T - \lambda I)u$$
\noindent for all $u \in D(T)$. Thus 
$$(T - \lambda I)Rg(T,\lambda)(T - \lambda I) = T - \lambda I$$
on $D(T)$.

On the other hand, if $u \in X$ and $v \in D(T)$ are 
such that $P(\lambda)u = (T - \lambda I)v$,
then
$$Rg(T,\lambda)(T - \lambda I)Rg(T,\lambda)u = 
 Rg(T,\lambda)(T - \lambda I)Q(\lambda)v$$
$$= Rg(T,\lambda)(T - \lambda I)v = Q(\lambda)v = Rg(T,\lambda)u.$$
\noindent Hence
$$Rg(T,\lambda)(T - \lambda I)Rg(T,\lambda) = Rg(T,\lambda)$$
\noindent on $X$.

Now we have to show that $Rg(T,\lambda) \in B(X)$. Because linearity is clear, and 
$D(Rg(T,\lambda)) = X$, it is sufficient to show that the operator $Rg(T,\lambda)$ is closed. 
Let $u_{n} \to u$ and
$Rg(T,\lambda)u_{n} \to w$ as $n \to +\infty$. Then there exist $v_{n} \in D(T)$ such that
$$P(\lambda)u_{n} = (T - \lambda I)v_{n} = (T - \lambda I)Q(\lambda)v_{n} \to P(\lambda)u.$$
On the other hand,
$$Rg(T,\lambda)u_{n} = Q(\lambda)v_{n} \to w = Q(\lambda)w,$$
since $R(Q(\lambda))$ is closed. Thus $Q(\lambda)v_{n} \to w$ and 
$(T - \lambda I)Q(\lambda)v_{n} \to P(\lambda)u$. Since $T$ is closed, we obtain that
$w \in D(T)$ and $P(\lambda)u = (T - \lambda I)w$. Hence $Rg(T,\lambda)u = Q(\lambda)w = w$,
showing that the operator $Rg(T,\lambda)$ is closed.

It is not complicated to see, using the definition of $Rg(T,\lambda)$, that
$$(T - \lambda I)Rg(T,\lambda) = P(\lambda)$$
and
$$Rg(T,\lambda)(T - \lambda I) = Q(\lambda)$$
for all $\lambda \in U$. We now show that $Rg(T,\lambda)$ satisfies the resolvent identity
in $U$. For $\lambda , \mu \in U$, we have
$$(\lambda - \mu)Rg(T,\lambda)Rg(T,\mu) = Rg(T,\lambda)\left[ \left( T - \mu I
 \right) - \left( T - \lambda I \right) \right] Rg(T,\mu)$$
$$= Rg(T,\lambda)P(\mu) - Q(\lambda)Rg(T,\mu)$$
$$= Rg(T,\lambda)P(\lambda)P(\mu) - Q(\lambda)Q(\mu)Rg(T,\mu)$$
$$= Rg(T,\lambda)P(\lambda) - Q(\mu)Rg(T,\mu)$$
$$= Rg(T,\lambda) - Rg(T,\mu).$$
\noindent Therefore $Rg(T,\lambda)$ verifies the resolvent identity in $U$ and (ii) is proved.

(ii) $\Rightarrow$ (i) : Take $P(\lambda) = (T - \lambda I)Rg(T,\lambda)$
and $Q(\lambda) = Rg(T,\lambda) (T - \lambda I)$ for all $\lambda \in U$.

(ii) $\Leftrightarrow$ (iii) : This equivalence follows from an easy
computation and using the fact that 
generalized resolvents are analytic \cite{HiPh}. \hfill $\diamondsuit$

\smallskip

\noindent {\bf REMARK 2.8} The projections $P(\lambda)$ and $Q(\lambda)$ obtained in (i) of the previous theorem
are analytic in $U$. Indeed, we have 
$P(\lambda) = (T - \lambda I)Rg(T,\lambda)$ and $Q(\lambda) = Rg(T,\lambda) 
(T - \lambda I)$. Since $Rg(T,\lambda)$ satisfies the resolvent identity, 
$P(\lambda)$ and $Q(\lambda)$ are 
analytic in $U$.


\section{The maximal radius of regularity of a Fredholm operator}
Let $T \in C(X)$. The operator $T$ is said to be a {\it Fredholm\/} operator
if $R(T)$ is closed and $\max$\{dim $N(T)$, codim $R(T)$\} $<$ $\infty$.
The Fredholm domain $\rho_{e}(T)$ of $T$ is defined by
$$\rho_{e}(T) = \{\lambda \in {\bf C} : T - \lambda I \   \mbox{Fredholm} \  \} .$$
The  set $\rho_{e}^{r}(T)$ defined by 
$$\rho_{e}^{r}(T) = \rho_{e}(T) \cap reg(T)$$ 
is an open set. 
Using \cite[Corollaire 2.3]{Mbe3} , we have that $\rho_{e}^{r}(T)$ is the 
set of all
$\lambda \in \rho_{e}(T)$ such that the application~$z \to$~dim~$N(T - zI)$ is 
constant in a neighborhood of $\lambda$. 
Using the continuity of the index, the same set $\rho_{e}^{r}(T)$ coincides with the set of 
all $\lambda \in \rho_{e}(T)$ such that the application~$z \to$~codim~$R(T - zI)$ 
is constant in a neighborhood of $\lambda$. Therefore
$$0 \in \rho_{e}^{r}(T) \Longrightarrow
dist(0, {\bf C}\setminus \rho_{e}^{r}(T)) = dist(0, \sigma_{g}(T)).$$
It is this distance that we call the maximal radius of regularity of $T$ (if 
$0 \in \rho_{e}^{r}(T)$).   

Let $\gamma (T)$ be the reduced minimum modulus
of $T$ : 
$$\gamma (T) = \inf \{ \| Tx \| : x \in D(T) \ ; \ dist(x,N(T)) = 1\}.$$
It was proved in 1975 by K.H. F\"{o}rster and M.A. Kaashoek \cite{FK} that 
$$(*) \quad 0 \in \rho_{e}(T) \Longrightarrow  dist(0, \sigma_{g}(T) \setminus \{ 0 \}) 
= \lim_{n \to \infty} 
\gamma (T^n)^{1/n}$$
and
$$(**) \quad 0 \in \rho_{e}^{r}(T) \Longrightarrow dist(0, \sigma_{g}(T)) 
= \lim_{n \to \infty} 
\gamma (T^n)^{1/n} .$$

Recently, J. Zem\'anek \cite{Zem} conjectured a different representation for 
the distance of $0$ to the left spectrum of $T$ if $T$ is 
assumed to be left invertible. Instead of the reduced minimum modulus, 
his representation is now in terms of the spectral radius of left inverses of 
$T$. Note that the conjecture in \cite{Zem} is stated in the 
more general framework of 
Banach algebras. To be more specific, let $\rho_{\ell}(T)$ be the 
set of all $\lambda \in {\bf C}$ such that $T - \lambda I$ is left invertible 
and let $\sigma_{\ell}(T) = {\bf C} \setminus \rho_{\ell}(T)$ be the 
left spectrum of $T$.  The {\it conjecture \/} in 
\cite{Zem} for bounded linear operators $T \in B(X)$ can be stated as follows : 
$$(***) \quad 0 \in \rho_{\ell}(T) \Longrightarrow dist(0,\sigma_{\ell}(T)) = 
\sup \{ \frac{1}{r(S)} : ST = I\},$$
where 
$r(S)$ is the spectral radius of $S$. 

One-sided invertible operators are particular cases of 
operators with generalized inverses. Therefore, we can consider the analogue 
of the conjecture of Zem\'anek for $dist(0,\sigma_{g}(T))$. The main results we will 
prove here are 
two analogues of the F\"{o}rster-Kaashoek's results for Zem\'anek's conjecture. 
This implies (Corollary \ref{cor2}) a positive answer for $(***)$ under the 
additional assumption of the Fredholmness of $T$.

We start with an alternative version of $(**)$.   

\begin{theorem} \label{the1}
Let $T \in C(X)$ be a linear operator such that $0 \in \rho_{e}^{r}(T)$.
Then 
$$ dist(0,\sigma_{g}(T)) = \sup \{ \frac{1}{r(S)} : TST = T\},$$
where $r(S)$ is the spectral radius of $S$.
\end{theorem}  

{\it PROOF. \/} Let $S \in B(X)$ be a generalized inverse of
$T$. Then, using Theorem \ref{cha1}, we get $T^{n}S^{n}T^{n} = T^{n}$, for all
$n \geq1$. By \cite[Lemma 4]{FK}, we have $\gamma (T^{n}) \geq 1/(\| S^{n} \|)$ and
therefore
$$\lim_{n \to\infty} \gamma (T^{n})^{1/n} \geq \frac{1}{r(S)}.$$
Using \cite[Theorem 5]{FK} we get the inequality
$$dist(0,\sigma_{g}(T)) \geq \sup\{ \frac{1}{r(S)} : TST = T\} .$$

 In order to prove the other inequality, set $d = dist(0,\sigma_{g}(T))$.
Then
$$B(0,d) = \{ \lambda \in {\bf C} : |\lambda | < d\} \subseteq \rho_{e}^{r}(T).$$

Suppose that $d < \infty $. Let $\varepsilon$ be a positive number and put
$$K = \overline{B}(0,\frac{d}{1+\varepsilon}) = \{\lambda \in {\bf C} : 
|\lambda | \leq \frac{d}{1+\varepsilon} \} .$$
Then $K \subset \rho_{e}^{r}(T)$ is compact. Using \cite[Theorem 3.1]{Mbe3},
there is a generalized resolvent $Rg(T,\lambda)$ for $T$ on $K$. Then,
for all $\lambda \in K$,
$$Rg(T,\lambda) = \sum_{n=0}^{\infty}\lambda ^{n}T_{n} , T_{n} \in B(X) , R(T_{n}) \subset D(T).$$
By Cauchy's integral formula, we have
$$T_{n} = \frac{1}{2\pi i}\int^{ }_{|\lambda | = \frac{d}{1+\varepsilon}}\lambda^{-(n+1)}Rg(T,\lambda)\,d\lambda$$
for all $n \geq 0$. Denoting $M = max\{ \| Rg(T,\lambda)\| : \lambda \in K\}$,
we obtain
$$\| T_{n}\| \leq M\left( \frac{1+\varepsilon}{d}\right)^{n+1}, n\geq 0.$$
Using now Theorem 2.6, we get
$$\| T_{0}^{n+1}\| \leq M\left( \frac{1+\varepsilon}{d}\right)^{n+1}, n\geq 0,$$
which implies $r(T_{0}) \leq (1+\varepsilon)/d$. Since $TT_{0}T = T$,
we have
$$\sup \{\frac{1}{r(S)} : TST = T\} \geq \frac{1}{r (T_{0})}
\geq \frac{d}{1+\varepsilon}.$$
Since $\varepsilon > 0$ was arbitrarily chosen, we have
$$\sup \{\frac{1}{r(S)} : TST = T\} \geq d .$$

Suppose now that $d = \infty$. Let $\varepsilon$ be a positive number and
$K = \{\lambda \in {\bf C} : |\lambda | \leq \frac{1}{\varepsilon}\}$. Then
$K \subset \rho_{e}^{r}(T)$ is compact. Using similar arguments and
notation, $T$ possess in $K$ a generalized resolvent $Rg(T,\lambda)$
and $r(T_{0}) \leq \varepsilon$. Hence
$$\sup \{\frac{1}{r(S)} : TST = T\} \geq \frac{1}{r(T_{0})}
 \geq \frac{1}{\varepsilon}.$$
Since this holds for every $\varepsilon > 0$, we have
$$\sup \{\frac{1}{r(S)} : TST = T\} = \infty = d.$$
The proof is now complete.  \hfill $\diamondsuit$

\smallskip


We recall now the Kato decomposition for Fredholm operators 
(cf. for instance \cite{FK}). Namely, $X$ decomposes into two $T$-invariant 
closed subspaces $X_0$ and $X_1$ and, if $T_i$ is the restriction of $T$ to 
$X_i$, $i = 0,1$, then :
\begin{enumerate}
\item dim $X_1 < \infty ;$ 
\item $T = T_0 \oplus T_1 ;$
\item $T_0$ is regular ;
\item $X_1 \subset N(T^k)$ for some $k$.
\end{enumerate}

The following result is a counterpart of $(*)$.   

\begin{theorem} \label{the2}
Let $T \in C(X)$ be a linear operator such that $0 \in \rho_{e}(T)$ and 
$X = X_0 \oplus X_1$ a fixed Kato decomposition with 
respect to $T$.
Then 
$$ dist(0,\sigma_{g}(T) \setminus \{ 0 \}) 
= \sup \{ \frac{1}{r(S)} : TST = T ; SX_0 \subseteq X_0 \} .$$
\end{theorem}  

{\it PROOF. \/} Since $X_1$ in Kato decomposition is finite dimensional, 
$T -\lambda I$ is Fredholm if and 
only if $T_0 - \lambda I_0$ is Fredholm, where $\lambda \in {\bf C}$ and 
$I_0$ is the identity operator on $X_0$. But $0 \in \rho_e(T)$; thus 
$T_0$ is Fredholm. Using $(*)$ and \cite[page 125]{FK}, we have  
$$dist(0,\sigma_{g}(T) \setminus \{ 0 \}) = \lim_{n \to \infty} 
\gamma (T^n)^{1/n} = \lim_{n \to \infty} 
\gamma (T^n_0)^{1/n} = dist(0,\sigma_{g}(T_0) \setminus \{ 0 \}).$$
Since $T_0$ is Fredholm, the previous theorem and the previous equations 
yield
$$dist(0,\sigma_{g}(T) \setminus \{ 0 \}) = 
\sup \{ \frac{1}{r(S_0)} : T_0S_0T_0 = T_0\}.$$
Let $d = \sup \{ \frac{1}{r(S_0)} : T_0S_0T_0 = T_0\}$ and 
suppose that $d < \infty$. Let $\varepsilon$ be a positive number. Then there 
exists $S_0 = S_0(\varepsilon) \in B(X_0)$ such that $T_0S_0T_0 = T_0$ and 
$$d - \varepsilon < \frac{1}{r(S_0)} \leq d .$$
Since $X_1 \subset N(T^k)$ for some $k$, 
the operator $T_1$ is defined on all of the finite dimensional space $X_1$ and 
is nilpotent. Let $T_1 = PT_1^{'}P^{-1}$ with 
$$T^{'}_1 = \left(\begin{array}{cccc}
			0 & 1 & \cdots & 0 \\
			\vdots & \ddots & \ddots & \vdots \\
			0 & \cdots  & 0 & 1
		\end{array}\right) \ .$$
Then  
$$S^{'}_1 = \left(\begin{array}{cccc}
			0 &   0 & \cdots & 0 \\
   1 & 0 & \cdots   & 0		\\
			\vdots & \ddots  & \ddots & \vdots \\
			0 & \cdots				  & 1 & 0
		\end{array}\right) \  $$
is a generalized inverse for $T^{'}_1$, $S_1 = PS^{'}_1P^{-1}$ is a 
generalized inverse for $T_1$ and $S_1$ is nilpotent. Let $S = S_0 \oplus 
S_1$ with respect to Kato's decomposition. Then $TST = T$ and 
$SX_0 \subseteq X_0$. In particular, the set in the $\sup$ is always nonvoid.

Since 
$S^n = S_0^n \oplus 
S_1^n$,for all $n$, and $S_1$ is nilpotent, we have $r(S) = r(S_0)$. Then 
$$d - \varepsilon < \frac{1}{r(S_0)} = \frac{1}{r(S)} \leq 
\sup \{ \frac{1}{r(S)} : TST = T ; SX_0 \subseteq X_0\}.$$
But this holds for all $\varepsilon > 0$ and thus 
$$d \leq \sup \{ \frac{1}{r(S)} : TST = T ; SX_0 \subseteq X_0 \}.$$

For the other inequality, consider an operator $S$ such that $TST = T$ and 
$SX_0 \subseteq X_0$. Then the matrix of $S$, with respect to Kato decomposition, 
has the following form

$$\left(\begin{array}{cc}
			\displaystyle{S_1} & \ast \\
			0 & \ast
		\end{array}\right) \ $$ 

and thus

$$ S^n = \left(\begin{array}{cc}
			\displaystyle{S_1^n} & \ast \\
			0 & \ast
		\end{array}\right) \ .$$

This implies $\| S^n \| \geq \| S_1^n \|$, so $r(S) \geq r(S_1)$. On the other 
hand, $TST = T$ implies $T_0S_1T_0 = T_0$. Therefore
$$ \frac{1}{r(S)} \leq \frac{1}{r(S_1)} \leq 
\sup \{ \frac{1}{r(S_0)} : T_0S_0T_0 = T_0 \} .$$
Thus $d \geq \sup \{ \frac{1}{r(S)} : TST = T ; SX_0 \subseteq X_0 \}.$ 

The case $d = \infty$ can be proved in a 
similar fashion.   \hfill $\diamondsuit$
\smallskip

\noindent {\bf REMARK} We do not know if condition $SX_0 \subseteq X_0$ can be removed 
in the above theorem. In follows from the same Theorem that the $\sup$ does not 
depend upon Kato decomposition.

\smallskip

Set $s(T) = \sup \{\frac{1}{r(S)} : TST = T\}$.

\begin{cor} \label{cor1}
Let $T \in C(X)$ and suppose that $0 \in \rho_{e}^{r}(T)$. Then, for all
$n \geq 1$, we have $s(T^{n}) = s(T)^{n}$.
\end{cor}

{\it PROOF. \/} Let $n \geq 1$.
Using Theorem 3.1 and \cite[Theorem 5]{FK}, 
we get
$$s(T) = \lim_{k\to \infty} \gamma (T^{k})^{1/k}.$$
Therefore $$s(T^{n}) = \lim_{k\to \infty} (\gamma (T^{kn})^{1/kn})^{n} 
= s(T)^{n}.$$ 
The proof is complete. \hfill $\diamondsuit$
\smallskip

We want to note that C. Schmoeger \cite{Sch} studied 
							when $\lim_{n\to \infty}s(T^{n})^{1/n}$ can be expressed as a 
distance
							from $0$ to a modified spectrum. Corollary \ref{cor1} implies his 
result in our
							situation.

The following consequence of Theorem \ref{the2} is a partial result for  
Zem\'anek's conjecture.

\begin{cor} \label{cor2}
Suppose that $T \in B(X)$ is a left invertible and Fredholm bounded 
linear operator. Then 
$$dist(0,\sigma_{\ell}(T)) = \sup \{ \frac{1}{r(S)} : ST = I\}.$$
\end{cor}

{\it PROOF. \/} We have $\sigma_g(T) \subset \sigma_{\ell}(T)$ and 
$0 \in \rho_{\ell}(T) \cap \rho_e(T)$. Therefore
$$dist (0, \sigma_{\ell}(T)) \leq dist (0,\sigma_{g}(T)) = 
\sup \{ \frac{1}{r(S)} : TST = T\}$$
\begin{center}
(by Theorem \ref{the2}) 
\end{center} 
$$ = \sup \{ \frac{1}{r(S)} : ST = I\}$$
\begin{center}
(since $T$ is left invertible)
\end{center}
$$ \leq dist (0, \sigma_{\ell}(T))$$
\begin{center}
(cf. \cite{Zem}).
\end{center}
The proof is complete. \hfill $\diamondsuit$

\smallskip

Now we mention and briefly discuss some open problems.

\smallskip

\noindent {\bf PROBLEM 3.5} When is the $\sup$ attained in the formula 
$$\sup \{\frac{1}{r(S)} : TST = T\} = 
dist(0,\sigma_{g}(T))$$
in Theorem \ref{the1} ? 

\smallskip

The proof of Theorem \ref{the1} shows that 
the $\sup$ is attained if and only if one can 
construct a generalized resolvent $Rg(T,\lambda)$ for $T$ defined on 
the {\sl all} 
open set $\rho_{e}^{r}(T)$, instead of arbitary compact subsets $K$. 
The later is 
still an open problem.

Note also that the formula of Theorem \ref{the1} 
implies the following result : If $T \in C(X)$, $0 \in \rho_{e}^r(T)$ and 
there exists a generalized inverse $S_0$ of $T$ such that $TS_0T = T$ and 
$\sigma(S_0) = \{ 0 \}$, then $\rho_{e}^{r}(T) = {\bf C}$. 
{\sl We think that the converse also 
holds.} This certainly holds if the answer to Problem 3.5 is positive. 
\smallskip

\noindent {\bf EXAMPLE 3.6} Consider $X = C[0,1]$ and $T \in C(X)$ given by 
$T(f) = f'$ and $D(T) = C^{1}[0,1]$. Then 
$\rho^{r}_{e}(T) = {\bf C}$. Moreover, there exists a generalized resolvent on the  whole 
${\bf C}$. Indeed, if we define 
$$Rg(\lambda)(f)(x) = \int_0^x f(t)e^{\lambda (x - t)}\, dt, \; \; f \in X ,$$ 
then we have 
$$P(\lambda) = (T - \lambda)Rg(\lambda) = I$$
and
$$Q(\lambda) = Rg(\lambda)(T - \lambda) = I - F(\lambda), $$
where $F(\lambda)(f)(x) = f(0)\exp(\lambda x)$. It is easy to see that both 
conditions of Theorem \ref{cha2} are satisfied. Thus  
$Rg(\lambda)$ satisfies the resolvent identity on ${\bf C}$. Also, the Volterra 
operator $S_0$ defined 
by 
$$S_0(f)(x) = Rg(0)(f)(x) = \int_0^x f(t)\, dt$$ 
is a quasi-nilpotent operator. \hfill $\diamondsuit$
\smallskip

Note that operators whose Fredholm domain is the whole complex plane 
were studied by several 
authors : cf. \cite{KaLa} and the references cited therein. 
We also want to note that the equality $\rho_{e}^{r}(T) = {\bf C}$ clearly implies that 
$\rho_{e}(T) = {\bf C}$. In fact, if $0 \in \rho_{e}(T)$, and $S_0 \in B(X)$ is a 
generalized inverse of $T$, then $\rho_{e}(T) = {\bf C}$ if and only if 
$\sigma_{e}(S_0) = \{0\}$, that is $S_0$ is a Riesz operator. Here $\sigma_{e}$ is the essential spectrum. 
Indeed, we have 
$$
(I - \lambda S_0)T = (T - \lambda S_0 T) = T - \lambda + \lambda (I - S_0 T) \; .
$$ 
Since $\lambda (I - S_0 T)$ is of finite rank, $T - \lambda$ is Fredholm 
if and only if $I - \lambda S_0$ is (cf. \cite{Gol}).
Therefore $\rho_{e}(T) = {\bf C}$ if and only if $S_0$ is a Riesz operator.

\smallskip

\noindent {\bf PROBLEM 3.7} If $0 \in reg(T)$, does it follow that 
							$s(T) = dist(0,\sigma_{g}(T))$ ?

\smallskip

For Hilbert space bounded operators the answer is positive. 
The details will be published elsewhere.

\vskip 1truecm


\bigskip


\noindent C. Badea \hfill M. Mbekhta 

\noindent URA D751 au CNRS \& \hfill  URA D751 au CNRS \& 

\noindent UFR de Math\'ematiques \hfill UFR de Math\'ematiques

\noindent Universit\'e de Lille I \hfill  Universit\'e de Lille I

\noindent F--59655 Villeneuve d'Ascq, France \hfill  F--59655 Villeneuve d'Ascq, France

\smallskip

\noindent e-mail address : \hfill e-mail address : 

\noindent {\tt badea@gat.univ-lille1.fr}  \hfill  {\tt mbekhta@gat.univ-lille1.fr}

\smallskip 

\noindent   $\; $ \hfill and

\smallskip

\noindent  $\; $ \hfill Universit\'e de Galatasaray, 

\noindent $\; $ \hfill \c{C}iragan Cad no 102

\noindent  $\; $ \hfill Ortakoy 80840, Istanbul, Turkey

\bigskip

\noindent 1991 Mathematical Subject Classification : 47A10, 47A05, 47A53.

\end{document}